\documentclass[11pt]{article}

\usepackage{amsmath}
\usepackage{amssymb}
\usepackage{amscd}
\usepackage{amsthm}

\newtheorem{theorem}{Theorem}[section]
\newtheorem{lemma}[theorem]{Lemma}
\newtheorem{proposition}[theorem]{Proposition}

\newtheorem{claim}[theorem]{Claim}
\newtheorem{definition}[theorem]{Definition\rm}
\newenvironment{example}{\textit{Example}}{\newline}
\newenvironment{examples}{\textit{Examples}}{\newline}
\newcommand{\A}{{\mathfrak A}}
\newcommand{\im}{{\operatorname{Im}}}
\begin{document}
\date{}
\title{A Representation Theorem for Completely Contractive
Dual Banach Algebras}

\author{Faruk Uygul}

\maketitle

\begin{abstract}
In this paper, we prove that every completely contractive dual
Banach algebra is completely isometric to a $w^\ast$-closed
subalgebra the operator space of completely bounded linear operators
on some reflexive operator space.
\end{abstract}

\section{Introduction}

A Banach algebra $\A$ which is dual Banach space is called a dual
Banach algebra if the multiplication on $\A$ is separately
$w^*-$continuous. All von Neumann algebras are dual Banach algebras,
but so are all measure algebras $M(G)$, where $G$ is a locally
compact group, and all algebras $\mathcal{B}(E)$, where $E$ is a
reflexive Banach space. Of course, every $w^\ast$-closed subalgebra
of $\mathcal{B}(E)$ for a reflexive Banach space $E$ is then also a
dual Banach algebra. Surprisingly, as proven recently by Daws
(\cite{Daw}), every dual Banach algebra arises in this fashion.

A completely contractive dual Banach algebra is an algebra which is
a dual operator space in the sense of \cite{E-R} such that
multiplication is completely contractive and separately
$w^\ast$-continuous. Then von Neumann algebras are examples of
completely contractive dual Banach algebras. Also, whenever $\A$ is
a dual Banach algebra, then $\max \A$ (\cite{E-R}) is a completely
contractive dual Banach algebra. If $G$ is a locally compact group,
then the Fourier--Stieltjes algebra $B(G)$ (\cite{Eym}) is a
completely contractive dual Banach algebra which, for arbitrary $G$,
is neither a von Neumann algebra nor of the form $\max B(G)$. In the
present paper, we prove an operator space analog of Daws'
representation theorem: if $\A$ is a completely contractive dual
Banach algebra, then there is a reflexive operator space $E$ and a
$w^\ast$-$w^\ast$-continuous, completely isometric algebra
homomorphism from $\A$ to ${\mathcal CB}(E)$, where ${\mathcal
CB}(E)$ stands for the completely bounded operators on $E$ (see
\cite{E-R}). We would like to stress that even where $\A$ is of the
form $\max \A$ for some dual Banach algebra $\A$, our result is not
just a straightforward consequence of Daws' result, but requires a
careful adaptation of his techniques to the operator space setting.
The construction of such a reflexive operator space heavily relies
on the theory of real and complex interpolation of operator spaces
defined by Xu (\cite{Xu}) and Pisier (\cite{Pis 1} and \cite{Pis 2})
respectively.

This representation theorem is somewhat related in spirit to results
by Ghahramani (\cite{Ghah}) and Neufang, Ruan, and Spronk
(\cite{N-R-S}): In \cite{Ghah}, $M(G)$ is (completely) isometrically
represented on $\mathcal{B}(L^2(G))$, and in \cite{N-R-S}, a similar
representation is constructed for the completely contractive dual
Banach algebra $M_{cb}(A(G))$. We would like to emphasize, however,
that our representation theorem neither implies nor is implied by
those results: $\mathcal{B}(L^2(G))$ is a dual operator space, but
not reflexive.

This paper is part of the author's PhD thesis written at the
University of Alberta under the supervision of Volker Runde. While
being in the final stages of the preparation of this manuscript, we
learned that our representation theorem had also been independently
proven by Daws.

\section{Preliminaries}
\subsection{Dual Banach Algebras and Operator Spaces}
\begin{definition}
A Banach algebra $\A$ is called a {\it dual Banach algebra\/} if it
is a dual Banach space and the multiplication on $\A$ is separately
$w^\ast$-continuous.
\end{definition}
\begin{example}1.\ Every von Neumann algebra is a dual Banach
algebra. \\ 2.\ If $E$ is a reflexive Banach space, then
$\mathcal{B}(E)$ is a dual Banach algebra with the predual $E^{\ast}
\hat{\otimes} E$, where $\hat{\otimes}$ represents the projective
tensor
product of Banach spaces. \\
3.\ If $G$ is a locally compact group, then the measure algebra
$M(G)$ and the Fourier-Stieltjes algebra $B(G)$ are dual Banach
algebras with preduals $C_0(G)$ and $C^\ast(G)$ respectively.
\end{example}
\begin{definition}
An {\it operator space\/} is a linear space $E$ with a complete norm
$\| \cdot \|_n$ on $M_n(E)$ for each $n \in \mathbb{N}$ such that
\begin{equation} \tag{R 1} \label{R1}
\left\| \begin{array}{c|c} x & 0 \\ \hline 0 & y \end{array}
\right\|_{n+m} = \max \{ \| x \|_n, \| y \|_m \} \qquad (n,m \in
{\mathbb N}, \, x \in M_n(E), \, y \in M_m(E))
\end{equation}
and
\begin{equation} \tag{R 2} \label{R2}
\| \alpha  x  \beta \|_n \leq \| \alpha \| \| x \|_n \| \beta \|
\qquad (n \in \mathbb{N}, \, x \in M_n(E), \, \alpha, \beta \in
M_n).
\end{equation}
\end{definition}
\begin{examples} 1.\ Every closed subspace of $\mathcal{B}(\mathcal{H})$,
where $\mathcal{H}$ is a Hilbert space, is an operator space.
\\ 2.\ If $G$ is a locally compact group, then
the measure algebra $M(G)$, the Fourier algebra $A(G)$, and the
Fourier-Stieltjes algebra $B(G)$ are operator spaces.
\end{examples}
\begin{definition}
Let $E_1,E_2$, and $F$ be operator spaces. A bilinear map $T \!: E_1
\times E_2 \to F$ is called {\it completely contractive\/} if \[ \|T
\|_{cb} := \sup_{n_1, n_2 \in \mathbb{N}} \left\|  T^{(n_1, n_2)}
\right\| \leq 1,\] where \\ $T^{(n_1, n_2)} : M_{n_1}(E_1) \times
M_{n_2}(E_2) \to M_{n_1n_2}(F), \quad  \big( (x_{i,j}), (y_{k,l})
\big) \mapsto \big( T(x_{i,j},y_{k,l})   \big)$.
\end{definition}
\begin{definition}
A {\it completely contractive Banach algebra\/} is an algebra which
is also an operator space such that multiplication is a completely
bounded bilinear map.
\end{definition}
\begin{examples} 1.\ Every closed subspace of $\mathcal{B}(\mathcal{H})$,
where $\mathcal{H}$ is a Hilbert space, is a completely contractive
Banach algebra.\\ 2.\ If $G$ is a locally compact group, then the
measure algebra $M(G)$, the Fourier algebra $A(G)$ and the
Fourier-Stieltjes algebra $B(G)$ are completely contractive Banach
algebras.
\end{examples}
\begin{definition}
Let $E$ and $F$ be operator spaces, and let $T \in
\mathcal{B}(E,F)$. Then:
\begin{enumerate}
\item $T$ is {\it completely bounded\/} if
\[
\| T \|_{\mathrm{cb}} := \sup_{n \in \mathbb{N}} \left\| T^{(n)}
\right\|_{{\mathcal B}(M_n(E),M_n(F))} < \infty.
\] \newline
\item $T$ is a {\it complete contraction\/} if $\| T \|_{\mathrm{cb}} \leq
1$.\newline
\item $T$ is a {\it complete isometry\/} if $T^{(n)}$ is an isometry for
each $n \in \mathbb{N}$. \newline
\end{enumerate}
The set of completely bounded operators from $E$ to $F$ is denoted
by $\mathcal{CB}(E,F)$.
\end{definition}
\subsection{Complex Interpolation of Banach Spaces}
Let $X_0, X_1$ be  two complex Banach spaces. The couple $(X_0,X_1)$
is called \textit{compatible} (in the sense of interpolation theory)
if there is a Hausdorff complex topological vector space
$\mathcal{X}$ and $\mathbb{C}$-linear continuous inclusions $X_0
\hookrightarrow \mathcal{X}$ and $X_1 \hookrightarrow \mathcal{X}$.

Now let $X_0$ and $X_1$ be two compatible normed spaces. Then we
define a norm on the set $X_0+X_1$ by $\|x\|:=\inf \{\|x_0 \|_{X_0}
+ \|x_1 \|_{X_1}: x=x_0+x_1 \}$. We denote this space by
$X_0+_1X_1$. For $ 0 < \theta <1 $, let
\[X_{[\theta]}=(X_0,X_1)_\theta:=\big\{x\in X_0+_1X_1: \ x=f(\theta) \ \rm{
\ for \ some \ f \
satisfying} \ (*) \ \big\} \] where \\
\textbf{* : }\ $f: \mathbb{C} \to X_0+_1X_1$ is a function which
satisfies the following conditions:\\
\\(i)\  $f$ is bounded and continuous on the strip $S:=\{z\in
\mathbb{C}: \ 0 \leq Re (z) \leq 1 \}$, \\
\\(ii)\ $f$ is analytic on $S_0$, the interior of $S$,\\
\\(iii) \  $f(it)\in X_0$ and $f(1+it)\in X_1 \ \ (t\in
\mathbb{R}).$\\
\\Define a norm on $X_{[\theta]}$ via
\[ \big\| x\big \|_{[\theta]}:= \inf
\Big\{\big \|f \big\| : \  x=f(\theta),\ f \ \rm{satisfying \ (*) }\
\Big\} \] where the norm of $f$ is defined to be
\[
\big\| f\big\|:= \max \Big\{\ \sup \big\{ \big\|f(it)\big \|_{X_0}
\big\}, \ \sup \big\{ \big\|f(1+it) \big\|_{X_1}\big\}: \ t\in
\mathbb{R} \Big\}.
\]
By this construction, $X_{[\theta]}$ becomes an interpolation space
between $X_0$ and $X_1$. For more information on interpolation of
Banach spaces, we refer the reader to \cite{B-L}.

\subsection{Interpolation of Operator Spaces}
Suppose that $E_0$, $E_1$ are operator spaces such that  $(E_0,E_1)$
is a compatible couple of Banach spaces. Note that for each $n \in
\mathbb{N}$, we have continuous inclusions $M_n(E_0) \hookrightarrow
M_n(\mathcal{\mathcal{X}})$ and $M_n(E_1) \hookrightarrow
M_n(\mathcal{\mathcal{X}})$ where $M_n(\mathcal{X})$ is identified
with $\mathcal{X}^{n^2}$. Thus $(M_n(E_0),M_n(E_1))$ is a compatible
couple of Banach spaces. Clearly $E_0 \oplus _{\infty} E_1$ is an
operator space by setting
\[
M_n(E_0 \oplus _{\infty} E_1)= M_n(E_0) \oplus _{\infty} M_n(E_1).
\]
Then $E_0 \oplus _{1} E_1$ becomes an operator space with the
embedding $E_0 \oplus _{1} E_1 \hookrightarrow (E_0^* \oplus
_{\infty} E_1^*)^*.$ Now, for each $1 < p < \infty$,  $E_0 \oplus
_{p} E_1$ becomes an operator spaces via
\[
M_n(E_0 \oplus _{p} E_1)= (M_n(E_0 \oplus _{1} E_1) ,M_n(E_0 \oplus
_{\infty} E_1))_\theta, \quad \quad 1/p=1- \theta .
\]

The complex and real interpolation of operator spaces are defined by
G. Pisier in \cite{Pis 1}   and by Q. Xu in \cite{Xu} respectively.
The construction of the latter heavily uses the first one.
\subsubsection{Complex Interpolation of Operator Spaces}
Let $0 < \theta < 1$ and $E_0$, $E_1$ be operator spaces such that
$(E_0,E_1)$ is a compatible couple of Banach spaces. Then, for each
$n \in \mathbb{N}$, the couple $(M_n(E_0), M_n(E_1))$ is also
compatible. Now define
\begin{equation} \label{complex interpolation}
M_n(E_\theta) := ( M_n(E_0), M_n(E_1))_\theta
\end{equation}
in the sense of complex interpolation (\cite{B-L}). By this
definition, $E_\theta = (E_0, E_1)_\theta$ becomes an operator
space. This is called the complex interpolation of operator spaces
$E_0$ and $E_1$ (see  \cite{Pis 1} and \cite{Pis 2} for more
information).

\subsubsection{Real Interpolation of Operator Spaces}
The construction of the interpolation of operator spaces by the real
method is more complicated than by the complex method. This is
because definition (\ref{complex interpolation}) does not work for
the real interpolation $(E_0,E_1)_{\theta , p}$ \ if $p < \infty$.
Now we introduce real interpolation of operator spaces by the
discrete $K$-method as defined by Xu in \cite{Xu}.

Note that if $E$ is an operator space and $t>0$, then $tE$ denotes
the operator space obtained by multiplying the norm on each matrix
level by $t$. Now let $\mu$ denote a weighted counting measure on
$\mathbb{Z}$ (That is: Let $\{a_n\}_{n \in \mathbb{Z}}$ be a
sequence of non-negative reals. For $E \subseteq \mathbb{Z} $, we
define $\mu (E) := \sum_{n \in E} a_n$ )  and $\{ E_k \}_{k \in
\mathbb{Z}}$ a sequence of operator spaces. Then for $1 \leq p \leq
\infty$, we define
\[
l_p (\{ E_k \}_{k \in \mathbb{Z}} ; \mu) := \{ (x_k)_{k \in
\mathbb{Z}}: \ \  x_k \in E_k \ \ \text{and} \ \ (\| x_k\|)_{k \in
\mathbb{Z}} \in l_p (\mu)      \}.
\]
Clearly $l_\infty (\{ E_k \}_{k \in \mathbb{Z}} ; \mu)$ is an
operator space with its natural operator space structure. Then $l_1
(\{ E_k \}_{k \in \mathbb{Z}} ; \mu)$ becomes an operator space when
it is considered as a subspace of $(l_\infty (\{ E^*_k \}_{k \in
\mathbb{Z}} ; \mu))^*.$ Finally $l_p (\{ E_k \}_{k \in \mathbb{Z}} ;
\mu)$ becomes an operator space by complex interpolation:
\[
l_p (\{ E_k \}_{k \in \mathbb{Z}} ; \mu) = (l_1 (\{ E_k \}_{k \in
\mathbb{Z}} ; \mu), l_\infty (\{ E_k \}_{k \in \mathbb{Z}} ; \mu)
)_\theta, \quad 1/p =1-\theta .
\]
Let $(E_0,E_1)$ be a compatible couple of operator spaces. For $1
\leq p \leq \infty$, we define  $N_p(E_0,E_1):=\{(x,-x): \ x\in E_0
\cap E_1 \}$ regarded as a subspace of $E_0 \oplus _p E_1$. Then we
define
\[
E_0 +_p E_1 := (E_0 \oplus _p E_1) / N_p(E_0,E_1) .
\]
$K_p(t;E_0,E_1)$ denotes the operator space $E_0 +_p tE_1 ;$ for any
$x\in E_0+E_1$, we let $K_p(x,t;E_0,E_1):= \|  x\|_{E_0 +_p tE_1 }$.
Now we may give the definition of $E_{\theta,p;\underline{K}}$, the
real interpolation of the compatible couple $(E_0,E_1)$ with the
discrete $K$-method, as follows:
\[
E_{\theta,p;\underline{K}} =(E_0,E_1)_{\theta,p;\underline{K}} =
\Big \{ x\in E_0+E_1: \ \ \|x\|_{\theta,p;\underline{K}} := \Big
[\sum_{k\in \mathbb{Z}}(2^{-k\theta}K_p(x,2^k;E_0,E_1))^p\Big]^{1/p}
\ < \infty \Big\}.
\] Then $E_{\theta,p;\underline{K}}$ is a Banach space.\\
If  $\alpha \in \mathbb{R}$, then $l_p(2^{k\alpha})$ is the weighted
space
\[
l_p(2^{k\alpha}):= \Big\{x=(x_k)_{k \in \mathbb{Z}}: \ \ \|
x\|_{l_p(2^{k\alpha})} = \Big(\sum_{k \in \mathbb{Z}} | 2^{k\alpha}
x_k |^p \Big)^{1/p} < \infty \Big\}.
\]
If $E$ is an operator space, then we similarly define $l_p(E)$ and
$l_p(E;2^{k\alpha})$ of sequences with values in $E$. Then $l_p(E)$
and $l_p(E;2^{k\alpha})$ are operator spaces.\\

For each $k \in \mathbb{Z}$, let  $F_k:=K_p(2^k;E_0,E_1).$ Then we
define $E_{\theta,p;\underline{K}}$, {\it the operator space
interpolation of the couple $(E_0,E_1)$ by the discrete $K$-method,
\/} as a subspace of $l_p (\{ F_k \}_{k \in \mathbb{Z}} ;
2^{-k\theta})$ consisting of the constant sequences. More
explicitly, let $x=(x_{i,j})\in M_n((E_0,E_1)_{\theta,p
;\underline{K}})$ for some $n \in \mathbb{N}$. Then
\begin{eqnarray*}
\|x\|_{M_n((E_0,E_1)_{\theta,p ;\underline{K}})}:=\inf \Big \{
\big\| (u,v)\big\|_{M_n\big(l_p(E_0;2^{-k\theta})\oplus _p
l_p(E_1;2^{k(1-\theta)})\big )} : \ (u,v) \ \text{satisfying} \ (**)
\Big \}.
\end{eqnarray*}
$$
(**) \quad u=(u_{i,j})\in M_n\big(l_p(E_0;2^{-k\theta})\big), \
v=(v_{i,j})\in M_n\big(l_p(E_1;2^{k(1-\theta)})\big)$$ \ \text{where
\ each}$$  u_{i,j}=\big(u_{i,j}^k\big)_{k \in \mathbb{Z}} \ \
\text{and} \ \ v_{i,j}=\big(v_{i,j}^k\big)_{k \in \mathbb{Z}}  \
$$ \text{such that} $$ \ x_{i,j}= u_{i,j}^k+v_{i,j}^k, \ \
\text{for each} \ \ k\in \mathbb{Z},\ \ i,j=1, \ldots n .
$$

\section{Main Theorem}
\textit{Notation.} \ 1. Let $(E_\alpha)_{\alpha \in I}$ be a family
of operator spaces where $I$ is some index set. Then $l^2-\bigoplus
_\alpha E_\alpha$ \ and\ $l^2(I,E_\alpha)$ will denote the $l^2$-
direct sum of $E_\alpha$'s and the complex operator space
interpolation $( l^\infty(I,E_\phi), l^1(I,E_\phi) )_{1/2}$
respectively. \\ 2. Let $\A$ be a completely contractive Banach
algebra and $X$ be an operator (bi- or) left $\A-$module. For
$a=(a_{i,j})$ and $x=(x_{i,j})$  in $M_n(\A)$ and $M_m(X)$
respectively, for some $n,m \in \mathbb{N},$ \ $x \star y$ will
represent the matrix
\begin{equation} \label{schur pro}
x \star y = (x_{i,j}.y_{k,l})
\end{equation}
where ``.'' represents the module action of $\A$ on $X$.
\begin{definition}\label{def2}
Let $\A$ be a completely contractive dual Banach algebra, $\phi \in
M_n(\A_*)$ for some $n\geq 1$. Suppose that for each $m\geq1,$ there
is a matricial norm $\|.\|_{\phi,m}$ on \ $M_m(\A .\phi)$. Let
$E_\phi$ denote the completion of \ $(\A . \phi , \| .
\|_{\phi,1}).$ Suppose that
\begin{eqnarray}
\|a\star b\|_{\phi, mk} &\leq& \| a \|_m \| b \|_{\phi ,k}
\label{condition1} \\
\text{and} \nonumber \\ \| a \star \phi \|_{mn} &\leq& \| a
\star\phi \|_{\phi ,m} \leq \| a \|_m \| \phi\|_{n}
\label{condition2}
\end{eqnarray}
for all $a\in M_m(\A) $ and $b\in M_k(E_\phi),  m,k \in \mathbb{N}.$\\
Furthermore, suppose that $E_\phi$ is reflexive and the inclusion \
$\iota _\phi: E_\phi \rightarrow M_n(\A_*)$  injective. Then $(\|.
\|_{\phi,m})_{m=1}^{\infty}$ is called an {\it admissible operator
norm \/}for $\phi.$\\
\end{definition}
 Note that in the previous definition, the inequality
(\ref{condition1}) means that $E_\phi$ is an operator left
$\A$-module.\\
\\
\begin{example} Let $G_d$ be a locally compact discrete group and
$\phi=(\phi_{i,j})\in M_n(C^*(G_d))$ with $ \|\phi\|_n=1 $ for some
$n \geq 1,$ \ where $C^*(G_d)$ is the group $C^*$-algebra of $G_d$.
Then each $\phi_{i,j}$ is a finite sum of the form
\[
\phi_{i,j} = \sum_{g\in G_d} \lambda ^{i,j}_g \delta _g,
\]
where each $\lambda^{i,j}_g \in \mathbb{C}$ and $\delta _g$ is the
Dirac function. Consider $E_\phi$ with the usual norm on
$M_n(C^*(G_d))$. Clearly $E_\phi$ is a closed subspace of
$M_n(C^*(G_d))$. Hence, it is an operator space. Clearly $E_\phi$ is
reflexive. Since $E_\phi$ is an operator $\A$-module,
(\ref{condition1}) is satisfied. Since $B(G_d)$ is a completely
contractive Banach algebra, (\ref{condition2}) is also satisfied.
Therefore, the usual norm on $M_n(C^*(G_d))$ defines an admissible
operator norm for $\phi.$
\end{example}

\begin{theorem} \label{proposition1}
Let $\A$ be a completely contractive  dual Banach algebra and let
$\phi =(\phi_{i,j}) \in M_n(\A_*)$ have an admissible operator norm
for some $n\geq 1$. Then there is a \ $w^*$-continuous, completely
contractive representation of $\A$ on $\mathcal{CB}(E_\phi).$
\end{theorem}
\begin{proof} It is easy to see that $E_\phi$ is a left $\A$-module.
Moreover, $\iota_\phi ^*$ has a dense range if and only if \
$\iota_\phi ^{**}: E_\phi^{**} \to M_n(\A^*)$ is injective. Since
$E_\phi$ is reflexive, $\iota_\phi ^{**} = \iota_\phi$. Hence,
$\iota_\phi ^*$ has a dense range. Note that
\[
\iota_\phi ^* : T_n(\A) \to E_\phi^*
\]
where $T_n(\A)$ is as defined in \cite{E-R}. Now we define
\begin{equation*} \label{S}
S_\phi: E_{\phi}^* \hat{\otimes} E_\phi \rightarrow M_{n^2}(\A_*),
\quad  \iota_\phi^*(b) \otimes a.\phi \mapsto a.\phi \star b,
\end{equation*}
where $\hat{\otimes}$ represents the projective tensor product of
operator spaces.\\
\\Due to Definition \ref{def2}, this map is well-defined. Then the
map defined by
\[
T_\phi:=S_\phi ^* : \ T_{n^2}(\A) \rightarrow CB (E_\phi)
\]
is $w^*$-continuous. Since $\A$ is completely isometrically
isomorphic to a closed subspace of $T_{n^2}(\A)$ by the map
\begin{equation} \label{embedding}
\A \to T_{n^2}(\A), \quad a \mapsto (a_{i,j})
\end{equation}
where \[ a_{i,j}= \left\{
\begin{array}{ll}
  a, & \hbox{if} \ \ (i,j)=(1,1); \\
  0, & \hbox{otherwise,}
\end{array}
\right. \] $T_\phi$ induces a multiplicative representation from
$\A$ into $CB(E_\phi)$. For simplicity, we will denote this
representation again by $T_\phi$. In order to see that $T_\phi$
is multiplicative on $\A$: \\
Let  $a,b \in \A$, $c=(c_{i,j}) \in T_n(\A)$. Consider $a$ as an
element of $T_{n^2}(\A)$ via (\ref{embedding}). Then
\[
\big \langle \langle T_\phi (a),b.\phi \rangle, \iota _\phi ^*
(c)\big \rangle = \langle T_\phi(a), \iota _\phi ^* (c) \otimes
b.\phi \rangle = \big \langle a, S_\phi (\iota _\phi ^* (c) \otimes
b.\phi) \big \rangle = \langle a, b.\phi \star c \rangle =\langle a,
b.\phi _{1,1}.c_{1,1} \rangle.
\]
Hence $\langle T_\phi (a),b.\phi \rangle = (x_{i,j})\in M_n(\A _*)$
where \[ x_{i,j}= \left\{
\begin{array}{ll}
  ab.\phi _{1,1}, & \hbox{if} \ \ (i,j)=(1,1); \\
  0, & \hbox{otherwise.}
\end{array}
\right. \] Now let $a,b,d \in \A $. Then
\[
\Big \langle T_\phi(d).T_\phi(a), b.\phi \Big \rangle = \Big \langle
T_\phi(d), \langle T_\phi(a), b. \phi \rangle \Big \rangle =
\Big\langle T_\phi(d), (x_{i,j}) \Big\rangle = (y_{i,j})= \Big
\langle T_\phi(da),b.\phi \Big \rangle \in E_\phi
\]
where
\[ y_{i,j}= \left\{
\begin{array}{ll}
  dab.\phi _{1,1}, & \hbox{if} \ \ (i,j)=(1,1); \\
  0, & \hbox{otherwise.}
\end{array}
\right. \] Hence $T_\phi$ is multiplicative on $\A$.\\

By using Effros and Ruan (\cite{E-R}, Proposition 7.1.2), $S_\phi$
is a complete contraction if and only if the induced map  \
$\widetilde{S}_\phi \in B^2(E_{\phi}^* \times E_\phi , M_{n^2}(\A_*)
)$ is a complete contraction. Now
\begin{eqnarray*}
\| \widetilde{S}_\phi \|_{cb}  & = &  \sup  \Big\{ \big\|
\widetilde{S}_{\phi}^{(m,m)} (x,y) \big\|  : x=(x_{i,j})\in
M_m(E_{\phi}^* ),\ y=(y_{i,j})\in M_m(E_{\phi}), \  \| x\|_{\phi
,m}, \\ & & \quad \quad \| y\|_{\phi ,m } \leq 1 , \ m \in
\mathbb{N} \Big\}
\\ & = & \sup \Big\{\big |\langle \! \langle
\widetilde{S}_{\phi}^{(m,m)} (x,y), z \rangle \! \rangle \big |:
x=(x_{i,j})\in M_m(E_{\phi}^* ),\ y=(y_{i,j})\in M_m(E_{\phi}),\\ &
& \quad \quad z=(z_{i,j})\in M_{m^2}(M_{n^2}(\A)),\ \| x\|_{\phi
,m}, \ \| y\|_{\phi ,m} , \ \|z\| _{m^2n^2} \leq 1, \ m \in
\mathbb{N} \Big \}.
\end{eqnarray*}
By the density of the range of $\iota_\phi ^*,$  for each
$i,j=1,\ldots m$, without loss of generality we may suppose that
\[
x_{i,j}= \iota_\phi^* (A_{i,j}), \quad y_{i,j}= b_{i,j}.\phi
\]
where
\[
A_{i,j}=(a_{i,j}^{k,l})\in T_n(\A) \ \text{and} \ \ b_{i,j}\in \A.
\]
Then we have
\begin{eqnarray*}  \Big \langle \! \Big\langle \widetilde{S}_\phi^{(m,m)}
(x,y), z  \Big\rangle  \! \Big\rangle   = \Big (\Big\langle
\widetilde{S}_\phi(x_{i,j},y_{k,l}),z_{s,t} \Big\rangle\Big)=
\big(\big\langle b_{k,l}.\phi \star A_{i,j},z_{s,t}
\big\rangle\big).
\end{eqnarray*}
Since
\[
\big\langle b_{k,l}. \phi _{o,p}. a_{i,j}^{q,r},z_{s,t} \big\rangle
= \big\langle \phi _{o,p}. a_{i,j}^{q,r},z_{s,t} b_{k,l}\big \rangle
= \big\langle a_{i,j}^{q,r},z_{s,t} b_{k,l}.\phi _{o,p} \big\rangle
\]
for all indices \ $i,j,k,l,m,n,o,p,q$ and $r$ \ where
\[
i,j,k,l= 1,\ldots m,  \quad  o,p,q,r= 1,\ldots n, \quad \text{and}
\quad s,t= 1,\ldots m^2n^2 ,
\]
we conclude that
\[
\Big |\Big \langle \! \Big\langle \widetilde{S}_\phi^{(m,m)} (x,y),
z \Big\rangle \! \Big\rangle\Big |  =  \Big |\Big\langle
\iota_\phi^* (A_{i,j}),z_{m,n}b_{k,l}.\phi \Big\rangle\Big |  = \Big
|\Big\langle \! \Big\langle (\iota_\phi^* (A_{i,j})),z \star y \Big
\rangle \! \Big\rangle \Big |.
\]
On the other hand,
\[
\| z \star y \|_{\phi,m^3n^2} \leq \| z \|_{m^2n^2} \| y \|_{\phi,m}
\leq 1.
\]
Therefore,\[ \Big|\Big\langle \! \Big\langle \widetilde{S}^{m,m}
(x,y), z \Big\rangle \!\Big \rangle\Big | \leq \Big \|\Big(
\iota_\phi^*(A_{i,j})\Big )\Big \|_{\phi ,m} = \| x  \|_{\phi ,m}
\leq 1.\] Thus, $\widetilde{S}_\phi$ is a complete contraction.
\end{proof}
$$$$
Let $(E_\alpha)_{\alpha \in I}$ be a family of operator spaces and
let $E=l_2(I,E_\alpha).$ We will need an approximation for the norm
of an arbitrary element on each matrix level of $E$. To manage this,
we will need the following two propositions.
\\
\begin{proposition} \label{prop 2 parts}
If $(E_\alpha)_{\alpha \in I}$ is a family of operator spaces, then\\
\\1. \ $l^{\infty}(I,M_n(E_\alpha))\cong
M_n(l^{\infty}(I,E_\alpha))$ \ for every $n \in \mathbb{N}$.\\
\\2. \ If $A=(a_{i,j}) \in M_n(l^{1}(I,E_\alpha)) \ \text{where} \
a_{i,j}=(a_{i,j}^\alpha)_\alpha, \ a_{i,j}^\alpha \in E_\alpha \
\text{for each $i$ and $j$}.$ Then
\[
\|   A \|_{M_n(l^{1}(I,E_\alpha))} \leq \sum _{\mathcal \alpha} \|
(a_{i,j}^\alpha) \|_{M_n(E_\alpha)} \ \text{ for every} \ n \in
\mathbb{N}.
\]
\end{proposition}
\begin{proof}
The first identity is obvious. Hence, we will prove only the second
one. Let $A$ be as in the claim. Then we have
\begin{eqnarray*}
\big\|  A \big\|_{M_n(l^{1}(I,E_\alpha))} & = & \sup \Big\{
\big|\big\langle \! \big\langle A,F \big\rangle \! \big\rangle
\big|: \  F \in M_n \big(l^{\infty}(I,E_\alpha^*)\big), \| F\|_{M_n
(l^{\infty}(I,E_\alpha^*))} \leq 1 \Big\} \\ & = & \sup \Big\{
\big|\big(\langle a_{i,j},f_{k,l} \big\rangle)\big |: \ F \in M_n
\big(l^{\infty}(I,E_\alpha^*)\big), \| F\|_{M_n
(l^{\infty}(I,E_\alpha^*))} \leq 1\Big \} \\
& = & \sup \Big\{ \big|\big(\sum _\alpha \langle f_{k,l}^\alpha ,
a_{i,j}^\alpha \big\rangle\big) \big|: \ F \in M_n
\big(l^{\infty}(I,E_\alpha^*)\big), \| F\|_{M_n
(l^{\infty}(I,E_\alpha^*))} \leq 1\Big \} \\ & \leq & \sup \Big\{
\sum _\alpha\big|\big ( \langle f_{k,l}^\alpha , a_{i,j}^\alpha
\rangle \big)\big|: \ F \in M_n (l^{\infty}(I,E_\alpha^*)), \|
F\|_{M_n (l^{\infty}(I,E_\alpha^*))} \leq 1 \Big\}
\\ & \leq& \sum _\alpha
\big\| (a_{i,j}^\alpha) \big \|_{M_n(E_\alpha)}.
\end{eqnarray*}\end{proof}

Next proposition will be needed to prove Lemma \ref{lemma} and
Theorem \ref{th 2}.
\\
\begin{proposition}  \label{BS Interpolation}
Let $(X,Z)$ be a compatible couple of Banach spaces in the sense of
Banach space interpolation. Suppose that there is a contractive
embedding from a Banach space $Y$ into $Z$. Let $E_1:=(X,Y)_\theta$
and $E_2:=(X,Z)_\theta$ for some $0 < \theta <1$. Then for every
$a\in E_1$, we have  $\|a\|_{E_1}\geq \|a\|_{E_2}$.
\end{proposition}
\begin{proof}
Clearly, any function $f: \mathbb{C} \to X+_1Y$ satisfying (*) can
be viewed as a function from $\mathbb{C}$ to $X+_1Z$, and it will
satisfy analogous properties. To distinguish these two functions, we
will denote them by $f_{X,Y}$ and $f_{X,Z}$
respectively. \\
Recall that
\begin{eqnarray*}
\big\| f_{X,Y}\big\|&=& \max \Big\{\ \sup \big\{ \big\|f(it)\big
\|_X \big\}, \ \sup \big\{
\big\|f(1+it) \big\|_Y\big\}: \ t\in \mathbb{R} \Big\}  \ \ \\
\rm{and} \\
\big\| f_{X,Z}\big\|&=&\max \Big\{\ \sup \big\{ \big\|f(it)
\big\|_X\big\},\ \sup \big\{ \big\|f(1+it) \big\|_Z\big\}:\ t\in
\mathbb{R}\Big \}.
\end{eqnarray*}
This shows that $E_1\subset E_2.$ \\
Now let $a\in E_1.$ Recall that  \[ \big\| a\big \|_{E_1}:= \inf
\Big\{\big \|f _{X,Y}\big\|: a=f(\theta),\ f \ \rm{satisfying \ (*)
}\ \Big\}.
\] ($\| a \|_{E_2}$ is defined in a similar fashion). Since for each
$f$ satisfying $(*)$ we have, $\| f_{X,Y} \| \geq \| f_{X,Z} \|$, we
conclude that $\|a \|_{E_1} \geq \|a \|_{E_2}$.
\end{proof}
Proposition \ref{prop 2 parts} and Proposition \ref{BS
Interpolation} show that:
\begin{lemma} \label{lemma}
Let $(E_\alpha)_{\alpha \in I}$ be a family of operator spaces and
let $A=(a_{i,j})$ be in $M_n(E)$ where $E=l_2(I, E_\alpha),$ for
some $n \in \mathbb{N}$. If each $a_{i,j}=(a_{i,j}^\alpha)_\alpha
\in E$, then $\|A\|_{M_n(E)} \leq \sqrt{\sum _\alpha
\|A_\alpha\|_{M_n(E_\alpha)}^2}$ \text{where each} $A_\alpha =
(a_{i,j}^\alpha)\in M_n(E_\alpha).$ That is, the canonical inclusion
from $ l_2(I, M_n(E_\alpha))$ into $M_n(l_2(I, E_\alpha))$ is a
contraction.
\end{lemma}

\begin{theorem}
Let $\A$ be a completely contractive dual Banach algebra. Then for
each $n \in \mathbb{N},$ every non-zero element in the unit ball of
$M_n(\A_*)$ has an admissible operator norm.
\end{theorem}
\begin{proof}
Let $\A$ be a completely contractive dual Banach algebra and $\phi
\in M_n(\A_*), \ \phi \neq 0, \ \| \phi \|_n \leq 1, $ for some $n
\geq 1$. The map
\[R_\phi : \A \to \A.\phi, \quad a \mapsto a.\phi\] is a complete
contraction. Then the induced map $\pi : \A /\ker R_\phi \to A.\phi$
is a complete isometry. For each $m \geq 1$, define a norm $\|
.\|_{\A.\phi, m}$ on $M_m(\A . \phi)$ via
\[
\|x\|_{\A.\phi, m} := \inf \big\{ \|a\|_m : \ \ x= a \star \phi, \ a
\in M_m(\A) \big\}.
\]
Then $\A . \phi$ becomes an operator space with this matricial norm.\\
Clearly $(\A .\phi, M_n(\A_*)) $ is a compatible couple of operator
spaces. Now define $E_\phi$ to be the space of constant sequences in
$l_2(\{F_k\}_{k \in\mathbb{N}}; 2^{-k/2})$, where $F_k=K_2(2^k;
\A.\phi, M_n(\A _*))$ for each $k\in \mathbb{N}$. By (\cite{B},
Section 2.3, Proposition 1) we know that:\begin{eqnarray*} (\A. \phi
, M_n(\A _*))_{\frac{1}{2},2;\underline{K}} \ \text{is reflexive}
\iff \text{the
inclusion} \  \A .\phi \to M_n(\A_*) \text{ is weakly compact} \\
\iff \text{the map} \ R_\phi : \A \to M_n(\A ^*), \ \ a \mapsto
a.\phi \ \text{is weakly compact}.\end{eqnarray*} However, $\im
(R_\phi) \subseteq M_n(\A_*)$ and $M_n(\A _*) \subseteq WAP
(M_n(\A^*))$, by \cite{Run 2}. Hence $R_\phi$ is weakly compact.
Therefore, as a closed subspace of $E=(\A .\phi,
M_n(\A_*))_{\frac{1}{2},2;\underline{K}}$, $E_\phi$ is reflexive
too.
\\Let $\| . \|_{\phi , m}$ denote the norm on $M_m(E_\phi)$, for
every $m \in \mathbb{N}$. If  $f\in E_\phi$, then
\begin{eqnarray*}
\| f \|_{\phi , 1}  & = & \Big [ \sum_{k \in \mathbb{N}}
2^{-k}\|f\|_{F_k}^2 \Big ]^{1/2} = \Big[ \sum_{k \in \mathbb{N}}
2^{-k} \inf _{b\in \A}\big\{ \sqrt{\|b.\phi\|_{\A .\phi ,1}^2 +
2^{2k}\|f-b.\phi \|_{n}^2 } \big\}^2\Big ]^{1/2}
\\ & = &  \Big [\sum_{k \in \mathbb{N}} 2^{-k} \inf _{b\in \A}\big\{
\|b.\phi\|_{\A .\phi ,1}^2 + 2^{2k}\|f-b.\phi \|_{n}^2   \big \}\Big
]^{1/2}
\\ & = &  \Big [ \sum_{k \in \mathbb{N}}  \inf _{b\in \A}\big\{ 2^{-k}
\|b.\phi\|_{\A .\phi ,1}^2 + 2^{k}\|f-b.\phi \|_{n}^2 \big \}\Big ]
^{1/2}.
\end{eqnarray*}
Hence,
\[
f \in E_\phi \iff \sum_{k \in \mathbb{N}}  \inf \big \{ 2^{-k}
\|b.\phi\|_{\A .\phi ,1}^2 + 2^{k}\|f-b.\phi \|_{n}^2 : \ b\in \A
\big \}  \ < \infty .
\]
Thus there exists a sequence $(b_k)_k$ in $\A$ such that \ $2^{2k}\|
f-b_k.\phi \|_{n}^2 \to 0.$ Hence $\A. \phi$ is dense in $E_\phi$.
This shows that $\| a. \phi \|_n \leq \| a. \phi \|_{\phi ,1}$. Now
to prove (\ref{condition2}), we will use the following claim.
\begin{claim} \label{c}
If \ $\A . \phi$ is dense in $E_\phi$, then $M_m(\A .\phi)$ is dense
in $M_m(E_\phi)$ for every $m \geq 1$.
\end{claim}
\begin{proof}
Let $\epsilon >0$ and $F=(f_{i,j}) \in M_m(E_\phi)$ for some $m \geq
1$. Then for each $i,j=1,\ldots m$, there exists a sequence
$(b_{i,j}^k)_k$ in $\A$ \ such that \ $b_{i,j}^k .\phi \to
f_{i,j}^k$ \ in  $\|.\|_{n}$. Consider the sequence $(F_k)_{k}$ in
$M_m(\A .\phi)$ where each $F_k=(b_{i,j}^k)_k$. Then we have
\[
\|F-f_k\|_{mn} \leq \sum_{i,j=1}^m \| f_{i,j}-b_{i,j}^k .\phi \|_n
\to 0 \ \ \ \text{as} \ \ \ k \to \infty .
\]
\end{proof}
Hence, by Claim \ref{c},  we have $\| a \star \phi \|_{mn} \leq \| a
\star \phi \|_{\phi,m} $ for all $a\in M_m(\A)$ and $m
\in \mathbb{N}.$\\
\\Let $a \in M_m(\A)$ for some $m \geq 1$. Then by the definition of
$\|.\|_{\A.\phi,m}$, it is clear that $\| a \star \phi \|_{\A.\phi ,
m} \leq  \|a \|_m$. Since $\A$ is a completely contractive dual
Banach algebra, we also have \  $ \|a \star \phi \|_{mn} \leq
\|a\|_m \|\phi \|_n \leq \| a\|_m$. By Lemma \ref{lemma}, we have
\[
\| a \star \phi \|_{F_k} \leq \inf \Big \{\Big[  \big  \| b \star
\phi\big \|_{\A . \phi ,m}^2 +2^{2k}\big\| a \star \phi -b \star
\phi\big \|_{mn}^2
      \Big ]^{1/2} : \ b \in M_m(\A) \Big \}.
\]
By choosing $b=a$, we see that
\[
\| a \star \phi \|_{F_k} \leq \| a \star \phi \|_{\A .\phi ,m} \quad
\text{for each $k$}.
\]
Then we have
\[
\| a \star \phi \|_{\phi, m} \leq \Big [\sum_{k \in \mathbb{N}}
2^{-k}\big\| a \star \phi \big\|_{F_k}^2\Big ]^{1/2} \leq \| a \star
\phi \|_{\A .\phi ,m} \Big [ \sum_{k \in \mathbb{N}}2^{-k}\Big]
^{1/2} = \| a \star \phi \|_{\A. \phi, m} \leq \| a \|_m.
\]
$$$$
Now let $a \in M_m(\A), b \in M_t(\A)$ for some \ $m,t \geq 1$.
Since the map $\pi$ is a complete isometry, we have
\begin{eqnarray*}
\big\| b \star (a \star \phi) \big\|_{\A.\phi, mt}  & = & \big  \|
b\star a + M_{mt}(\ker R_\phi)\big\|= \inf \big\{\| b\star a
+x\|_{mt} : \ x \in M_{mt}(\ker R_\phi)       \big\} \\ & \leq &
\inf\big \{\| b\star a +b \star x\|_{mt}: \ x \in M_{m}(\ker R_\phi)
\big\}\\& \leq& \| b \|_t \inf \big\{\| a+x\|_m: \ x \in M_{m}(\ker
R_\phi) \big\}
\\ & = & \| b \|_t \| a+\ker R_\phi \| = \|b \|_t \|a\star \phi
\|_{\A.\phi,m}.
\end{eqnarray*}
This shows that $\A.\phi$ is an operator left $\A-$module. Since
$M_n(\A_*)$ is also an operator left $\A-$module, so is $E_\phi$.
Therefore,
\[
\| b \star d \|_{\phi, mt} \leq \|b \|_t \|d \|_{\phi, m}  \quad
\text{for every} \ \ d \in M_m(E_\phi), \ b \in M_m (\A),\ m,t \in
\mathbb{N}.
\]
\end{proof}

\begin{theorem} \label{th 2}
Let $\A$ be a completely contractive  dual Banach algebra such that
for each $n \in \mathbb{N},$ every non-zero element in the unit ball
$\mathfrak{I}$ of $M_n(\A_*)$ has an admissible operator norm. Then
there is a $w^\ast$-continuous complete isometry map from $\A$ into
$CB(E)$ for some reflexive operator space $E$.\\
\end{theorem}
\begin{proof} Let $\dot{E}:= l^2-\bigoplus_{\phi \in \mathfrak{I}}E_\phi$.
For $n \geq 1$, equip $M_n(\dot{E})$ by $\| A \|_{M_n(\dot{E})}=
\sqrt{\sum_\phi \| (a_{i,j}^\phi) \|_{\phi , n}^2}$ where
$A=(a_{i,j})\in M_n(\dot{E}),$ $a_{i,j}=(a_{i,j}^\phi)_\phi,$
$a_{i,j}^\phi \in
E_\phi .$ \ Note that $\dot{E}$ is not an operator space.\\
There is a natural map \[S: \A \to B(\dot{E}) \ \ \rm{ defined \ by}
\ \ \ \Big\langle S(a),(x_\phi)\Big\rangle := \Big(\langle T_\phi
(a), x_\phi \rangle\Big),\] where $T_\phi: \A \to CB(E_\phi)$ is the
$w^*$-continuous complete contraction as defined in Theorem
\ref{proposition1}. \\Note that Daws (\cite{Daw}, Theorem 4.5)
proved that $S: \A \to B(\dot{E})$
is an isometry.\\
For an arbitrary  $n \geq 1,$ any element $(a_{i,j})$ of  $M_n(\A)$
can be viewed as a map \[S_n: \dot{E} \to M_n(\dot{E}).\] We claim
that $S_n$ is a contraction. Let $\| (a_{i,j}) \|_n \leq 1$ and
$(x_\phi)\in \dot{E}$. Then
\begin{eqnarray*}
\|S_n(x_\phi)\|_{M_n(\dot{E})} & = &  \Big\| \big(\big\langle
S(a_{i,j}),(x_\phi)\big\rangle\big) \Big\|_{M_n(\dot{E})} =
\Big\|\big (\big(\big\langle T_\phi(a_{i,j}),x_\phi
\big\rangle\big)\big)\Big \|_{M_n(\dot{E})} \\ &=& \Big[\sum_\phi
\big\| ((\langle T_\phi(a_{i,j}),x_\phi \rangle))
\big\|_{\phi,n}^2\Big ]^{1/2}   \leq  \Big [ \sum_\phi
\big\|(a_{i,j}) \big\|_n^2 \big\|x_\phi \big\|_{\phi,1}^2\Big]^{1/2}
\\ &=& \big\| (a_{i,j}) \big\|_n \Big [\sum_\phi \| x_\phi
\|_{\phi,1}^2\Big]^{1/2}= \big\| (a_{i,j}) \big\|_n \big\| (x_\phi)
\big\|_{\dot{E}} \leq \big\| (x_\phi) \big\|_{\dot{E}}.
\end{eqnarray*}
Thus, $S_n$ is a contraction.
\\  Now let
$E=l^2(\mathfrak{I},E_\phi)$. We define \begin{equation} \label{def
T} T: \A \to CB(E) \quad \rm{by}\ \ \ \langle T(a),(x_\phi) \rangle
= (\langle  T_\phi(a),x_\phi \rangle).\end{equation} Then \
$T^{(n)}: M_n(\A) \to CB (E,M_n(E))$.  On the other hand, by its
definition \[ M_n(E)= ( M_n(l^\infty(\mathfrak{I},E_\phi)),
M_n(l^1(\mathfrak{I},E_\phi)) )_{1/2}  = (
l^\infty(\mathfrak{I},M_n(E_\phi)), M_n(l^1(\mathfrak{I},E_\phi))
)_{1/2}.\] Each $(a_{i,j})\in M_n(\A) $ defines a map from $E$ into
\[\Big( l^\infty\big(\mathfrak{I},M_n(E_\phi)\big),
l^1\big(\mathfrak{I},M_n(E_\phi)\big)\Big)_{1/2}= l^2-\bigoplus
_{\phi \in \mathfrak{I}} M_n(E_\phi) \quad ( \rm{on \ the \ Banach \
space \ level}). \] Hence, we have a natural map (which we will
denote by $\widetilde{T}^n$) \
\[\widetilde{T}^n : M_n(\A) \to B \big(E, l^2-\bigoplus_{\phi \in
\mathfrak{I}} M_n(E_\phi)\big).\] However, this map is a contraction
for every $n \geq 1$. On the other hand, by Proposition \ref{BS
Interpolation} we have $\|\widetilde{T}^n \| \geq \| T^{(n)}\|$.
Hence, $T^{(n)}$ is a contraction for every $n \geq 1$. Thus $T$ is
a complete contraction.

Note that without loss of generality we may suppose that $\A$ is
unital and let $e$ denote its identity. For each $n \geq 1$, we have
\[
T^{(n)}:  \ M_n(\A) \to M_n(CB(E))=M_n((E^* \hat{\otimes}
E)^*)=CB((E^* \hat{\otimes} E),M_n).
\]
Let $a=(a_{i,j})\in M_n(\A)$. Then for every $\epsilon > 0$, \ there
is \ $\phi = (\phi _{i,j}) \in M_n(\A_*)$ such that
\[
\| \phi \|_n \leq 1 \ \ \text{and} \ \    | \langle \! \langle
a,\phi \rangle \! \rangle | \geq (1-\epsilon) \| a \|_n.
\]
For simplicity, set $\overline{T} := T^{(n)}(a) \in CB((E^*
\hat{\otimes} E),M_n)$. Define $ x=(x_{i,j})\in M_n (E^*
\hat{\otimes} E) $  by
\[ x_{i,j}=
\left\{
\begin{array}{ll}
  (\ldots,\iota_\phi^* (B) ,\ldots) \otimes (\ldots,e.\phi,\ldots), &
\hbox{if} \ \ (i,j)=(1,1); \\
  0, & \hbox{otherwise,}
\end{array}
\right. \] where \[B=(b_{i,j})\in T_n(\A) \ \ \text{is\ defined\ by}
\ \ b_{i,j}= \left\{
\begin{array}{ll}
  e, & \hbox{if} \ \ (i,j)=(1,1); \\
  0, & \hbox{otherwise.}
\end{array}
\right. \] Now we have \[\|x \|_{M_n (E^* \hat{\otimes} E)}= \|x_1
\|_{E^* \hat{\otimes} E} \leq \| (\ldots,\iota_\phi^*(B) ,\ldots)
\|_{E^*} \| (\ldots,e.\phi,\ldots) \|_{E}.  \] On the other hand,
\[
\| (\ldots,\iota_\phi^* (B) ,\ldots) \|_{E^*} \leq 1 \ \ \text{
since} \ \  \| B \| _{T_n(\A)} \leq 1 \ \ \text{and} \ \
\iota_\phi^* \ \ \text{is a contraction.}
\]
Clearly
\[
\| (...,e.\phi,...) \|_E \leq \|\phi \|_n \leq 1.
\]
Then
\[
\big| \overline{T}^{(n)}(x)\big | =
\big|\big(\overline{T}(x_{1,1})\big)\big | =\big |\big(\langle
T(a_{i,j}),x_{k,l} \rangle\big)\big| = \big| \langle \! \langle a,
\phi \rangle \!\rangle\big | \geq (1-\epsilon) \| a\|_n.
\]
\end{proof}


\section{Acknowledgements}

I would like to thank my supervisor Volker Runde for his valuable
recommendations and giving me the reference to the paper
\cite{Xu}.\\

\smallskip { \sc Department of Mathematical and Statistical
Sciences,\\
University of Alberta,
Edmonton, Alberta Canada, T6G 2G1.}\\
\\ E-mail: {\tt fuygul@ualberta.ca}


\begin{thebibliography}{G--L--W}
%
\begin{small}
%


\bibitem [1] {B} B. BEAUZAMY,
{\it Espaces d'Interpolation R\'{e}els: Topologie et
G\'{e}om\'{e}trie.} Springer Verlag, 1978.

\bibitem [2] {B-L} J. BERGH and J. L\"{O}FSTR\"{O}M,
{\it Interpolation Spaces.} Springer Verlag, 1976.


\bibitem[3]{Daw} M. DAWS,
{\it Dual Banach Algebras: representations and injectivity.}  (to
appear).

\bibitem[4] {E-R} E.G. EFFROS and Z.-J. RUAN,
{\it Operator Spaces.} Oxford University Press, 2000.

\bibitem [5] {Eym} P. EYMARD,
{\it L'alg\`{e}bre de Fourier d'un groupe localement compact.} Bull.
Soc. Math. France {\bf 92} (1964), 181--236.


\bibitem[6] {Ghah} F. GHAHRAMANI,
{\it Isometric representations of M(G) on B(H).} Glasgow Math. J.,
\textbf{23} (1982), 119--122.


\bibitem[7] {N-R-S} M. NEUFANG, Z.-J. RUAN and N. SPRONK,
{\it Completely isometric representations of $M_ cb A (G)$ and
$mathit UCB (hatG)^*$.} Trans. Amer. Math. Soc. (to appear).




\bibitem[8] {Pis 1} G. PISIER,
{\it The operator Hilbert space OH, complex interpolation and tensor
norms.} Mem. Amer. Math. Soc., \textbf{585} (1996).

\bibitem[9] {Pis 2} G. PISIER,
{\it Non-commutative vector valued $L_p$-spaces and completely
$p$-summing maps.} Ast\'{e}risque., \textbf{247} (1998).



\bibitem[10] {Ruan} Z.-J. RUAN,
{\it Subspaces of $C^\ast-$algebras.} J. Funct. Anal., \textbf{76}
(1988), 217--230.


\bibitem[11] {Run} V. RUNDE,
{\it Amenability for dual Banach algebras.} Studia Math.,
\textbf{148} (2001), 47--66.

\bibitem[12] {Run 2} V. RUNDE,
{\it Dual Banach algebras: Connes-amenability, normal, virtual
diagonals, and injectivity of the predual bimodule.} Math. Scand.,
\textbf{95} (2004), 124--144.



\bibitem[13] {Xu} Q. XU,
{\it Interpolation of Operator Spaces.} J. Funct. Anal.,
\textbf{139} (1996), 500--539.


\end{small}
%
\end{thebibliography}
\end{document}